\documentclass[11pt]{article}
\usepackage[margin=1in]{geometry}
\usepackage{amsmath,amssymb,amsthm,mathtools,microtype}
\usepackage[T1]{fontenc}
\usepackage{lmodern}
\usepackage{xcolor}
\usepackage[colorlinks=true,hypertexnames=false,linkcolor=blue!50!black,urlcolor=blue!50!black,citecolor=blue!50!black]{hyperref}

\newtheorem{theorem}{Theorem}[section]
\newtheorem{lemma}[theorem]{Lemma}
\newtheorem{proposition}[theorem]{Proposition}
\theoremstyle{remark}
\newtheorem{remark}[theorem]{Remark}

\newcommand{\R}{\mathbb R}
\newcommand{\E}{\mathbb E}
\newcommand{\Pc}{\mathcal P_c(\R^d)}
\newcommand{\Pt}{\mathcal P_2(\R^d)}
\newcommand{\norm}[1]{\lVert #1\rVert}
\newcommand{\OT}{\operatorname{OT}_{\varepsilon}}
\newcommand{\KL}{D_{\mathrm{KL}}}
\DeclareMathOperator{\Cov}{Cov}
\DeclareMathOperator{\tr}{tr}
\DeclareMathOperator{\supp}{supp}
\title{Long-Time Error Bounds for the Sinkhorn Wasserstein Gradient Flow}
\author{Alice Xuming Feng\thanks{Department of Management Science and Engineering, Stanford University. Email: \href{mailto:alisfeng@stanford.edu}{\texttt{alisfeng@stanford.edu}}.}}
\date{{\normalfont\normalsize (This is a preliminary draft. Comments are welcome.)}}

\begin{document}
\maketitle

\begin{abstract}
We study the long-time behavior of the Wasserstein gradient flow of the Sinkhorn divergence on $\Pt$.
In the absence of general convergence guarantees, we investigate how closely the flow approximates the target law at large times.
For initial and target laws $q_0,p\in\Pc$ with $q_0\ll\mathcal L^d$, we establish the bounds $\limsup_{T\to\infty}T^{-1}\int_0^T W_2^2(q_t,p)\,dt\leq\varepsilon d$ and $\limsup_{t\to\infty}W_2^2(q_t,p)\leq C\varepsilon d$, where $\varepsilon$ is the entropic regularization parameter and $C<5.013$ is a universal constant.
In deriving these estimates, we establish a general inequality bounding the quadratic Wasserstein distance in terms of the Sinkhorn divergence for arbitrary probability measures in $\Pt$.
\end{abstract}

\section{Introduction}

Recent applications in generative modeling~\cite{wflow,zhu-neural} motivate the study of the Wasserstein gradient flow (WGF) of the Sinkhorn divergence $S_\varepsilon(\cdot,p)$, where $p$ is a fixed target law. A basic question is whether the evolving law $q_t$ converges to $p$ as $t\to\infty$.

One route to convergence is through geodesic convexity. Under the standard hypotheses of Ambrosio, Gigli, and Savar\'e~\cite{ags}, if an energy is $\lambda$-convex along generalized geodesics with $\lambda>0$ and has a minimizer $q^*$, its Wasserstein gradient flow satisfies $W_2(q_t,q^*)\le e^{-\lambda t}W_2(q_0,q^*)$. For the Sinkhorn divergence (as the energy functional), however, although $p$ is the unique minimizer of $S_\varepsilon(\cdot,p)$~\cite{feydy}, the available geodesic convexity estimate has a negative constant, rendering the above argument inapplicable~\cite{ccl}. For Gaussian initial and target laws, Hardion and Lacombe~\cite{hardion-lacombe} proved convergence of WGF of Sinkhorn divergence whenever the initial covariance is nonsingular. They also raised the possibility that an initial density might suffice beyond the Gaussian setting.

This paper studies the long-term behavior of Sinkhorn WGF of general compactly supported initial and target laws. We ask how far the flow can remain from the target at long times. With the regularization parameter $\varepsilon>0$ held fixed, can we bound this distance in terms of $\varepsilon$ and the dimension $d$ alone?

For compactly supported laws $q_0$ and $p$ on $\R^d$, with $q_0$ absolutely continuous, we obtain
\[
\limsup_{T\to\infty}\frac{1}{T}\int_0^T W_2^2(q_t,p)\,dt\le\varepsilon d,
\qquad
\limsup_{t\to\infty}W_2^2(q_t,p)<5.012112\,\varepsilon d.
\]

Our results show that the long-time Wasserstein error can be made arbitrarily small by choosing a sufficiently small regularization parameter, even without establishing exact convergence for any fixed \(\varepsilon>0\).
\subsection{Related work}

Ambrosio, Gigli, and Savar\'e~\cite{ags} developed a general theory of gradient flows in metric spaces and in the space of probability measures. Under their hypotheses, convexity along generalized geodesics yields existence, uniqueness, and stability of Wasserstein gradient flows. A positive convexity constant also gives exponential convergence to the minimizer.

For entropically regularized optimal transport, Carlier, Chizat, and Laborde~\cite{ccl} establish displacement smoothness and prove well-posedness of the Wasserstein gradient flow of the Sinkhorn divergence on compact convex domains. Han et al.~\cite{wflow} establish global well-posedness for compactly supported initial and target laws in $\R^d$, with support bounds on each finite time interval.

Hardion and Lacombe~\cite{hardion-lacombe} study convergence of the Sinkhorn Wasserstein gradient flow in the Gaussian setting. The closed-form expression for the entropy-regularized quadratic transport cost between Gaussian measures~\cite{mallasto-gaussian} allows them to characterize the flow through the evolution of its mean and covariance. They prove convergence for nonsingular initial Gaussian laws, exhibit failures of convergence for singular initial laws, and derive convergence rates when the covariance matrices commute. In Section~6 of their paper, they suggest that an initial density might suffice for convergence beyond the Gaussian setting. Our work bounds the long-time Wasserstein distance to the target for general compactly supported laws without requiring a parametric description of the flow.

In generative modeling, Genevay, Peyr\'e, and Cuturi~\cite{genevay-generative} use Sinkhorn divergences as losses for fitting generative models. Zhu et al.~\cite{zhu-neural} approximate the velocity field of the Sinkhorn Wasserstein gradient flow using neural networks, while Han et al.~\cite{wflow} employ Wasserstein gradient flows for one-step generation.

\section{Notation}

We write $\Pt$ for the probability measures on $\R^d$ with finite second moment, $\Pc$ for those with compact support, and $\mathcal L^d$ for Lebesgue measure. For a measure $\mu$ and a measurable set $A$, the restriction $\mu|_A$ is defined by $(\mu|_A)(B)=\mu(A\cap B)$ for every measurable set $B$. For a measurable map $X$, the pushforward $X_\#\mu$ is defined by $(X_\#\mu)(B)=\mu(X^{-1}(B))$.
The set of couplings of $\mu$ and $\nu$ is denoted by $\Pi(\mu,\nu)$, and the quadratic Wasserstein distance is
\[
W_2^2(\mu,\nu)=\inf_{\gamma\in\Pi(\mu,\nu)}\int\norm{x-y}^2\,d\gamma.
\]
Here $\norm{\cdot}$ is the Euclidean norm. Following Han et al.~\cite{wflow}, we write $B_R:=\{x\in\R^d:\norm{x}\le R\}$ for the closed ball of radius $R$ centered at the origin. On $L^2(\mu)$, we use the real inner product $\langle f,g\rangle_{L^2(\mu)}:=\int fg\,d\mu$.

For random variables $A,B$, their mutual information is $I(A;B):=\KL\bigl(\mathcal L(A,B)\Vert\mathcal L(A)\otimes\mathcal L(B)\bigr)$.
For a random vector $A$ with density $f$, its differential entropy is $h(A):=-\int f(x)\log f(x)\,dx$, whenever this integral is well defined.
For a joint law that depends on time, $\E_t$ and $\Cov_t$ denote expectation and covariance under the specified law at time $t$; the same convention applies to conditional expectations and covariances.
For example, if $(X,Y)\sim\pi^\varepsilon_{q_t,p}$, then $\E_t[f(X,Y)]=\int f(x,y)\,\pi^\varepsilon_{q_t,p}(dx,dy)$.

For a vector-valued map $F$, we write $DF$ for its Jacobian, with $(DF)_{ij}=\partial_{x_j}F_i$; subscripts on $D$ specify the differentiation variable when needed. Unless explicitly stated otherwise, $D$ denotes the classical derivative; weak derivatives will be indicated explicitly. We reserve $\nabla f$ for the gradient of a scalar-valued function and write $\nabla\cdot F=\tr(DF)$ for divergence.
For a linear map $A$ between Euclidean spaces, its operator norm is $\norm A_{\mathrm{op}}:=\sup_{\norm h=1}\norm{Ah}$.
This convention applies in particular to $DF(x)$.

For bounded self-adjoint operators $A,B$ on $L^2(\mu)$, we write $A\preceq B$ if $\langle f,Af\rangle_{L^2(\mu)}\le\langle f,Bf\rangle_{L^2(\mu)}$ for every $f\in L^2(\mu)$, and write $I$ for the identity operator. Operators act componentwise on vector-valued functions. In particular, for the coordinate map $\mathrm{id}(x)=x$, we use the shorthand $Ax=A\,\mathrm{id}$, so that $(Ax)_i(x)=[A(y\mapsto y_i)](x)$.

\section{Problem setup and main result}

Throughout, $p$ is the fixed target law, $q_t$ is the evolving law, and $\varepsilon>0$ is fixed along the flow.

For $\mu,\nu\in\Pt$, set
\begin{equation}
\OT(\mu,\nu):=
\inf_{\pi\in\Pi(\mu,\nu)}
\left\{\frac12\int\norm{x-y}^2\,d\pi
+\varepsilon\KL(\pi\Vert\mu\otimes\nu)\right\}.
\label{eq:ot}
\end{equation}
Here $\KL(\pi\Vert\alpha)=\int\log(d\pi/d\alpha)\,d\pi$ if $\pi\ll\alpha$, and $+\infty$ otherwise. Define the Sinkhorn divergence of Feydy et al.~\cite{feydy}
\[
S_\varepsilon(\mu,\nu)=\OT(\mu,\nu)
-\frac12\OT(\mu,\mu)-\frac12\OT(\nu,\nu).
\]
We denote the unique minimizer in~\eqref{eq:ot} by $\pi^\varepsilon_{\mu,\nu}$~\cite{peyre-cuturi}. For compactly supported $\mu,\nu$, let $(u,v)=(u^\varepsilon_{\mu,\nu},v^\varepsilon_{\mu,\nu})$ be a pair of optimal dual potentials. The optimal coupling satisfies
\begin{equation}
\frac{d\pi^\varepsilon_{\mu,\nu}}{d(\mu\otimes\nu)}(x,y)
=\exp\!\left(\frac{u^\varepsilon_{\mu,\nu}(x)+v^\varepsilon_{\mu,\nu}(y)-\norm{x-y}^2/2}{\varepsilon}\right).
\label{eq:gibbs-coupling}
\end{equation}
The pair $(u,v)$ is unique up to the shift $(u,v)\mapsto(u+c,v-c)$. We use the Schr\"odinger equations to extend the potentials smoothly to all of $\R^d$.
\[
\begin{aligned}
u^\varepsilon_{\mu,\nu}(x)
&=-\varepsilon\log\int\exp\!\left(\frac{v^\varepsilon_{\mu,\nu}(y)-\norm{x-y}^2/2}{\varepsilon}\right)\nu(dy),\\
v^\varepsilon_{\mu,\nu}(y)
&=-\varepsilon\log\int\exp\!\left(\frac{u^\varepsilon_{\mu,\nu}(x)-\norm{x-y}^2/2}{\varepsilon}\right)\mu(dx).
\end{aligned}
\]
Write $\pi^\varepsilon_{\mu,\nu}(dx,dy)=\mu(dx)\,\pi^\varepsilon_{\mu,\nu}(dy\mid x)$. The barycentric projection is
\[
T^\varepsilon_{\mu,\nu}(x)=\int y\,\pi^\varepsilon_{\mu,\nu}(dy\mid x).
\]
For compactly supported $\mu$ and $p$, Han et al.~\cite{wflow} give the Sinkhorn velocity as the negative spatial gradient of the first variation:
\[
V^\varepsilon_{\mu,p}(x)
:=-\nabla_x\frac{\delta S_\varepsilon(\cdot,p)}{\delta\mu}(\mu)(x)
=T^\varepsilon_{\mu,p}(x)-T^\varepsilon_{\mu,\mu}(x).
\]
The population flow solves
\begin{equation}
\partial_tq_t+\nabla_x\cdot\bigl(q_tV^\varepsilon_{q_t,p}\bigr)=0.
\label{eq:sinkhorn-flow}
\end{equation}

\begin{theorem}[Long-term bias]\label{thm:main}
Let $d\ge1$, $\varepsilon>0$, and assume
\[
q_0\ll\mathcal L^d,\qquad q_0,p\in\Pc.
\]
The Sinkhorn continuity equation~\eqref{eq:sinkhorn-flow} with initial law $q_0$ admits a unique global solution $(q_t)_{t\ge0}$ with bounded support on each finite time horizon. Then
\begin{equation}
\limsup_{T\to\infty}\frac1T\int_0^T W_2^2(q_t,p)\,dt
\le\varepsilon d.
\label{eq:cesaro}
\end{equation}
Furthermore,
\begin{equation}
\limsup_{t\to\infty}W_2^2(q_t,p)
<5.012112\,\varepsilon d.
\label{eq:main}
\end{equation}

\end{theorem}

\begin{remark}[Scope of the assumptions]\label{rem:scope}
The AC assumption is only needed for $q_0$ and $p$ need not have a density. Compact support is used to ensure the existence and the uniqueness of the sinkhorn WGF. To the best of our knowledge, well-posedness of the Sinkhorn Wasserstein gradient flow for general initial and target laws with unbounded support remains an open question~\cite[Section~1]{hardion-lacombe}. Theorem~\ref{thm:main} bounds the long-time Wasserstein error but does not establish convergence of $q_t$.
\end{remark}

\section{Well-posedness and energy dissipation}

We first justify the existence and uniqueness assertion in Theorem~\ref{thm:main} and establish the characteristic representation used in the subsequent estimates.

We first establish the following identity.

\begin{lemma}\label{lem:jac}
Let $\mu,\nu\in\Pc$ and $(X,Y)\sim\pi^\varepsilon_{\mu,\nu}$. Then, for every $x\in\R^d$, using the conditional law in~\eqref{eq:gibbs-conditional},
\begin{equation}
D_xT^\varepsilon_{\mu,\nu}(x)=\varepsilon^{-1}\Cov(Y\mid X=x).
\label{eq:jac-cov}
\end{equation}
\end{lemma}

\begin{proof}
For compactly supported $\mu,\nu$, the canonical conditional law is
\begin{equation}
\pi^\varepsilon_{\mu,\nu}(dy\mid x)
=\frac{\exp\!\left(\frac{v^\varepsilon_{\mu,\nu}(y)-\norm{x-y}^2/2}{\varepsilon}\right)\nu(dy)}
{\int\exp\!\left(\frac{v^\varepsilon_{\mu,\nu}(z)-\norm{x-z}^2/2}{\varepsilon}\right)\nu(dz)}
=\frac{e^{x\cdot y/\varepsilon}\,\beta(dy)}
{\int e^{x\cdot z/\varepsilon}\,\beta(dz)},
\label{eq:gibbs-conditional}
\end{equation}
where
\[
\beta(dy)=\exp\!\left(\frac{v^\varepsilon_{\mu,\nu}(y)-\norm y^2/2}{\varepsilon}\right)\nu(dy).
\]
The Gibbs formula~\eqref{eq:gibbs-conditional} gives
\[
(T^\varepsilon_{\mu,\nu})_i(x)
=\frac{\int y_i e^{x\cdot y/\varepsilon}\,\beta(dy)}
{\int e^{x\cdot y/\varepsilon}\,\beta(dy)}.
\]
Since $\beta$ is finite and compactly supported, the integrands and their $x_j$-derivatives are bounded by a $\beta$-integrable constant locally uniformly in $x$; the dominated convergence theorem therefore justifies differentiation under both integrals.
For $(X,Y)\sim\pi^\varepsilon_{\mu,\nu}$, the quotient rule gives
\[
\begin{aligned}
\partial_{x_j}(T^\varepsilon_{\mu,\nu})_i(x)
&=\frac1\varepsilon\left[
\frac{\int y_i y_j e^{x\cdot y/\varepsilon}\,\beta(dy)}
{\int e^{x\cdot y/\varepsilon}\,\beta(dy)}
-\frac{\left(\int y_i e^{x\cdot y/\varepsilon}\,\beta(dy)\right)
\left(\int y_j e^{x\cdot y/\varepsilon}\,\beta(dy)\right)}
{\left(\int e^{x\cdot y/\varepsilon}\,\beta(dy)\right)^2}
\right]\\
&=\varepsilon^{-1}\bigl(\E[Y_iY_j\mid X=x]
-\E[Y_i\mid X=x]\E[Y_j\mid X=x]\bigr)\\
&=\varepsilon^{-1}\Cov(Y_i,Y_j\mid X=x),
\end{aligned}
\]
which is the claimed matrix identity.
\end{proof}

\begin{theorem}[Well-posedness and characteristic representation]\label{thm:flow}
Assume $q_0,p\in\Pc$. Define the initial and target support radii by
\[
R_0:=\sup_{x\in\supp q_0}\norm{x},\qquad
R_p:=\sup_{y\in\supp p}\norm{y}.
\]
The continuity equation~\eqref{eq:sinkhorn-flow} admits a unique global solution $(q_t)_{t\ge0}$ with bounded support on each finite time horizon. For every $T>0$,
\begin{equation}
q_t=(X_t)_\#q_0,\qquad
\dot X_t(a)=V^\varepsilon_{q_t,p}(X_t(a)),\qquad X_0(a)=a,
\label{eq:characteristics}
\end{equation}
for $0\le t\le T$. Each $X_t$ is bi-Lipschitz on $\operatorname{int}B_{R_0+1}$. Moreover,
\[
\supp q_t\subseteq B_{R_T},\qquad
R_T=e^TR_0+R_p(e^T-1),\qquad 0\le t\le T.
\]
Moreover, if $q_0\ll\mathcal L^d$, then $q_t\ll\mathcal L^d$ for every finite $t$.
\end{theorem}

\begin{proof}
By Han et al.~\cite[Theorem~B.4 and Proposition~B.7]{wflow}, the continuity equation~\eqref{eq:sinkhorn-flow} has a unique global solution with bounded support on each finite horizon. Proposition~B.5 of the same paper gives the stated support bound.

Fix $T>0$. The barycentric formulas and the support bounds give
\[
\norm{V^\varepsilon_{q_t,p}(x)}\le R_p+R_T,
\qquad x\in\R^d,\quad t\in[0,T].
\]
By Lemma~\ref{lem:jac}, if $\supp\nu\subseteq B_R$, then for every unit vector $h$,
\[
h^\top D_xT^\varepsilon_{\mu,\nu}(x)h
=\varepsilon^{-1}\operatorname{Var}(h\cdot Y\mid X=x)
\le\varepsilon^{-1}\E[(h\cdot Y)^2\mid X=x]
\le R^2/\varepsilon.
\]
Since the covariance matrix is positive semidefinite, this bounds the operator norm of $D_xT^\varepsilon_{\mu,\nu}(x)$. Applying it to both barycentric terms gives
\[
\norm{D_xV^\varepsilon_{q_t,p}(x)}_{\mathrm{op}}
\le L_T:=\frac{R_p^2+R_T^2}{\varepsilon},
\qquad x\in\R^d,\quad t\in[0,T].
\]
Thus the velocity is globally $L_T$-Lipschitz in space.
For every closed ball $B$ containing $B_{R_T}\cup\supp p$, \cite[Lemma~B.6, equation~(40)]{wflow} gives
\[
\sup_{x\in B}\norm{V^\varepsilon_{q_t,p}(x)-V^\varepsilon_{q_s,p}(x)}
\le 3C_{B,\varepsilon}W_2(q_t,q_s)
\longrightarrow0\qquad(t\to s).
\]
Hence the velocity is continuous in time locally uniformly in space.
The global bounds above imply, for every $r>0$,
\begin{equation}
\int_0^T\left(
\norm{V^\varepsilon_{q_t,p}}_{L^\infty(B_r)}
+\mathop{\operatorname{ess\,sup}}_{x\in B_r}\norm{D_xV^\varepsilon_{q_t,p}(x)}_{\mathrm{op}}
\right)dt
\le T(R_p+R_T+L_T)<\infty,
\label{eq:velocity-local-sobolev}
\end{equation}
and
\begin{equation}
\int_0^T\mathop{\operatorname{ess\,sup}}_{x\in\R^d}
\frac{\norm{V^\varepsilon_{q_t,p}(x)}}{1+\norm{x}}\,dt
\le T(R_p+R_T)<\infty.
\label{eq:velocity-growth}
\end{equation}
These verify the local Sobolev and growth hypotheses of \cite[Remark~7]{ambrosio-crippa}. Therefore \cite[Proposition~4 and Remark~7]{ambrosio-crippa} apply and give the global characteristic flow and \eqref{eq:characteristics}.

For $a,b\in\R^d$, integrating the difference of the characteristic equations and using the global Lipschitz bound gives
\[
\norm{X_t(a)-X_t(b)}
\le\norm{a-b}+L_T\int_0^t\norm{X_s(a)-X_s(b)}\,ds.
\]
For fixed $t\in[0,T]$, integrating from $s$ to $t$ also gives
\[
\norm{X_s(a)-X_s(b)}
\le\norm{X_t(a)-X_t(b)}
+L_T\int_s^t\norm{X_u(a)-X_u(b)}\,du,
\qquad 0\le s\le t.
\]
Gronwall's inequality applied to the first inequality gives the upper bound below. For the lower bound, reverse time in the second inequality by setting $s=t-r$ and changing variables $u=t-v$ in the integral:
\[
\norm{X_{t-r}(a)-X_{t-r}(b)}
\le\norm{X_t(a)-X_t(b)}
+L_T\int_0^r\norm{X_{t-v}(a)-X_{t-v}(b)}\,dv,
\qquad 0\le r\le t.
\]
Applying Gronwall's inequality in the variable $r$ gives
\[
\norm{X_{t-r}(a)-X_{t-r}(b)}
\le e^{L_Tr}\norm{X_t(a)-X_t(b)}.
\]
Evaluating at $r=t$ and using $X_0(a)=a$ and $X_0(b)=b$ yields
$\norm{a-b}\le e^{L_Tt}\norm{X_t(a)-X_t(b)}$. Combining the two bounds, we obtain
\begin{equation}
e^{-L_Tt}\norm{a-b}\le\norm{X_t(a)-X_t(b)}
\le e^{L_Tt}\norm{a-b},\qquad a,b\in\R^d.
\label{eq:flow-bilipschitz}
\end{equation}
The lower bound implies that $X_t$ is injective. Therefore its inverse exists on its image and is Lipschitz, so it preserves Lebesgue-null sets. For any Borel set $A$ with $\mathcal L^d(A)=0$, this gives
\[
\mathcal L^d\!\left(X_t^{-1}(A)\right)=0,
\qquad
q_t(A)=q_0\!\left(X_t^{-1}(A)\right)=0
\]
whenever $q_0\ll\mathcal L^d$. Hence $q_t\ll\mathcal L^d$.

\end{proof}

For the rest of the paper, we abbreviate $V_t:=V^\varepsilon_{q_t,p}$.

\begin{lemma}[Energy dissipation]\label{lem:dissipation}
Along the flow of Theorem~\ref{thm:flow}, for every $T>0$,
\begin{equation}
S_\varepsilon(q_T,p)+\int_0^T\int\norm{V_t(x)}^2q_t(dx)\,dt
=S_\varepsilon(q_0,p).
\label{eq:edi}
\end{equation}
By nonnegativity of the Sinkhorn divergence, this gives
\begin{equation}
\int_0^\infty\int\norm{V_t(x)}^2q_t(dx)\,dt<\infty.
\label{eq:dissipation}
\end{equation}
\end{lemma}

\begin{proof}
Fix $T>0$, and choose a compact convex set $K$ such that
\[
\supp p\;\cup\!\bigcup_{0\le t\le T}\supp q_t\subset\operatorname{int}K,
\]
which is possible by Theorem~\ref{thm:flow}. Applying the two-sided estimate in \cite[proof of Theorem~4.1]{ccl} to the cross- and self-cost terms yields $\lambda_T$-convexity of $S_\varepsilon(\cdot,p)$ along generalized geodesics on $\mathcal P(K)$ for some $\lambda_T\in\R$. It is finite and continuous on $\mathcal P(K)$ and bounded below by zero. Extending it by $+\infty$ outside $\mathcal P(K)$ preserves lower semicontinuity. Thus the hypotheses of \cite[Theorem~11.2.1]{ags} hold.

The curve $q_t$ is $W_2$-absolutely continuous by the bounded velocity and characteristic representation. Indeed, the coupling $(X_s,X_t)_\#q_0$ gives, for $0\le s\le t\le T$,
\[
W_2(q_s,q_t)
\le\norm{X_t-X_s}_{L^2(q_0)}
\le\int_s^t\norm{V_u}_{L^2(q_u)}\,du
\le(\operatorname{diam}K)(t-s).
\]
Since its supports stay away from $\partial K$, the no-flux boundary condition imposed in \cite{ccl} on $K$ is automatically satisfied. Hence \cite[Proposition~4.3]{ccl} identifies $(q_t)_{0\le t\le T}$, with velocity $V_t$, as the Wasserstein gradient flow of $S_\varepsilon(\cdot,p)$ on $\mathcal P(K)$. The energy identity \cite[Theorem~11.2.1, equation~(11.2.4)]{ags} therefore gives
\[
S_\varepsilon(q_T,p)+\int_0^T\int\norm{V_t(x)}^2q_t(dx)\,dt
=S_\varepsilon(q_0,p),
\]
proving \eqref{eq:edi}.
Since $S_\varepsilon(q_T,p)\ge0$, \eqref{eq:edi} bounds the accumulated dissipation by $S_\varepsilon(q_0,p)$ for every $T>0$, and letting $T\to\infty$ proves \eqref{eq:dissipation}.
\end{proof}

\section{Additional auxiliary results}

\subsection{Auxiliary identities and self-plan estimates}\label{sec:self}

\begin{lemma}\label{lem:density}
Along the flow of Theorem~\ref{thm:flow}, let $r_0(dx)=\eta_0(x)\,dx$ be a probability law with $\supp r_0\subseteq\supp q_0$, and set $r_s=(X_s)_\#r_0$.
Then $r_s(dx)=\eta_s(x)\,dx$, where, for every finite $s\ge0$,
\begin{equation}
\eta_s(X_s(a))
=\eta_0(a)\exp\!\left(-\int_0^s\nabla_x\cdot V_u(X_u(a))\,du\right)
\qquad\text{for $r_0$-almost every $a$}.
\label{eq:liouville-restricted}
\end{equation}
Taking $r_0=q_0$ gives the density transformation for $q_s$ whenever $q_0\ll\mathcal L^d$.
\end{lemma}

\begin{proof}
Fix $T>0$. Estimates~\eqref{eq:velocity-local-sobolev} and~\eqref{eq:velocity-growth} verify the local Sobolev and growth hypotheses of \cite[Remark~7]{ambrosio-crippa} for $V$ on $[0,T]$. Since $X_t$ is bi-Lipschitz and $r_0\ll\mathcal L^d$, also $r_t=(X_t)_\#r_0\ll\mathcal L^d$; write $r_t=\eta_t\mathcal L^d$. By \cite[Proposition~4 and Remark~7]{ambrosio-crippa},
\[
\partial_t\eta+\nabla\cdot(\eta V)=0,
\qquad\text{hence}\qquad
\partial_t\eta+V\cdot\nabla\eta=-\eta\nabla\cdot V
\]
in distributions, with initial density $\eta_0$.
The bound $\norm{D_xV_t(x)}_{\mathrm{op}}\le L_T$ from the proof of Theorem~\ref{thm:flow} gives
\begin{equation}
|\nabla\cdot V_t(x)|
=|\tr(D_xV_t(x))|
\le \norm{D_xV_t(x)}_{\mathrm{op}}\,d
\le L_T\,d,
\qquad t\in[0,T],\quad x\in\R^d.
\label{eq:velocity-divergence-bound}
\end{equation}
Since $r_t$ is a probability measure, this implies
\[
\eta,\ \eta\nabla\cdot V\in L^1([0,T]\times\R^d)
\subset L^1_{\mathrm{loc}}([0,T]\times\R^d).
\]
Thus \cite[Proposition~6 and Remark~7]{ambrosio-crippa}, applied with $w=\eta$, $b=V$, and $c=-\eta\nabla\cdot V$, give a representative of $\eta$ such that, for $\mathcal L^d$-almost every $a$ and every $s\in[0,T]$,
\[
\eta_s(X_s(a))
=\eta_0(a)-\int_0^s\eta_u(X_u(a))\nabla\cdot V_u(X_u(a))\,du.
\]
The lower bound in \eqref{eq:flow-bilipschitz} gives $\operatorname{Lip}(X_u^{-1})\le e^{L_Tu}$. Thus, for every Borel set $A\subseteq\R^d$,
\[
((X_u)_\#\mathcal L^d)(A)
=\mathcal L^d(X_u^{-1}(A))
\le (e^{L_Tu})^d\mathcal L^d(A)
\le e^{L_TT\,d}\mathcal L^d(A).
\]
Hence $(X_u)_\#\mathcal L^d\le e^{L_TT\,d}\mathcal L^d$, which, using Tonelli's theorem to exchange the order of integration, gives
\[
\int_{\R^d}\int_0^T |\eta_u(X_u(a))\nabla\cdot V_u(X_u(a))|\,du\,da
\le e^{L_TT\,d}\int_0^T\int_{\R^d}|\eta_u(x)\nabla\cdot V_u(x)|\,dx\,du<\infty.
\]
This finite bound implies $u\mapsto\eta_u(X_u(a))\nabla\cdot V_u(X_u(a))\in L^1(0,T)$ for $\mathcal L^d$-almost every $a$. By the fundamental theorem of calculus for Lebesgue integrals, $s\mapsto\eta_s(X_s(a))$ is absolutely continuous for such $a$ and satisfies
\[
\frac{d}{ds}\eta_s(X_s(a))
=-\nabla\cdot V_s(X_s(a))\eta_s(X_s(a))
\quad\text{for almost every $s\in[0,T]$}.
\]
Solving this scalar linear ODE yields
\[
\eta_s(X_s(a))
=\eta_0(a)\exp\!\left(-\int_0^s\nabla\cdot V_u(X_u(a))\,du\right).
\]
Since $r_0\ll\mathcal L^d$ and $T$ is arbitrary, this proves \eqref{eq:liouville-restricted}.
\end{proof}

\begin{proposition}[Self-plan displacement and self-barycentric drift]\label{prop:self}
Let $\mu\in\Pc$, and let $(X,Z)\sim\pi^\varepsilon_{\mu,\mu}$. Then
\begin{equation}
\E\norm{X-Z}^2\le\varepsilon d.
\label{eq:self}
\end{equation}
\begin{equation}
\int\norm{x-T^\varepsilon_{\mu,\mu}(x)}^2\mu(dx)
\le\frac12\E\norm{X-Z}^2\le\frac{\varepsilon d}{2}.
\label{eq:self-drift}
\end{equation}
\end{proposition}

\begin{proof}
Because the cost and the two marginals are symmetric and the optimizer is unique, the self-entropic optimizer is symmetric. Choose the two Schr\"odinger potentials equal and define
\[
\alpha(dx):=\exp\!\left(\frac{u^\varepsilon_{\mu,\mu}(x)}{\varepsilon}\right)\mu(dx).
\]
Continuity of the potential and compactness of $\supp\mu$ ensure that $\alpha$ is a finite positive measure supported on $\supp\mu$. The Gibbs formula~\eqref{eq:gibbs-coupling} for the optimal coupling then gives
\begin{equation}
\pi^\varepsilon_{\mu,\mu}(dx,dz)
=e^{-\norm{x-z}^2/(2\varepsilon)}\alpha(dx)\alpha(dz).
\label{eq:schrodinger-factorization}
\end{equation}
For $s>0$, keep $\alpha$ fixed and define
\begin{equation}
F(s)=\iint e^{-s\norm{x-z}^2/2}\alpha(dx)\alpha(dz).
\label{eq:self-F-definition}
\end{equation}
At $s=\varepsilon^{-1}$, $F(\varepsilon^{-1})=1$ and
\[
-2F'(\varepsilon^{-1})=\E\norm{X-Z}^2.
\]
The Gaussian Fourier representation gives
\begin{equation}
e^{-s\norm{x-z}^2/2}
=(2\pi s)^{-d/2}\int_{\R^d}
e^{-\norm\xi^2/(2s)}e^{i\xi\cdot(x-z)}\,d\xi.
\label{eq:gaussian-fourier}
\end{equation}
Substituting \eqref{eq:gaussian-fourier} into \eqref{eq:self-F-definition}, the resulting integrand is absolutely integrable, since
\[
\begin{aligned}
&\int_{\R^d}\iint
\left|e^{-\norm\xi^2/(2s)}e^{i\xi\cdot(x-z)}\right|
\,\alpha(dx)\alpha(dz)\,d\xi\\
&\qquad=\alpha(\R^d)^2\int_{\R^d}e^{-\norm\xi^2/(2s)}\,d\xi<\infty.
\end{aligned}
\]
Fubini's theorem therefore permits exchanging the $x,z$-integrations with the $\xi$-integration, giving
\begin{equation}
\begin{aligned}
F(s)
&=(2\pi s)^{-d/2}\int_{\R^d}e^{-\norm\xi^2/(2s)}
\left(\int e^{i\xi\cdot x}\alpha(dx)\right)
\left(\int e^{-i\xi\cdot z}\alpha(dz)\right)d\xi\\
&=(2\pi s)^{-d/2}\int_{\R^d}
e^{-\norm\xi^2/(2s)}|\widehat\alpha(\xi)|^2\,d\xi.
\end{aligned}
\label{eq:self-fourier}
\end{equation}
Here $\widehat\alpha(\xi):=\int e^{i\xi\cdot x}\alpha(dx)$ is the Fourier transform of $\alpha$.
By \eqref{eq:self-fourier}, $s^{d/2}F(s)$ is nondecreasing in $s$. Differentiating therefore yields
\[
F'(s)+\frac d{2s}F(s)\ge0.
\]
At $s=\varepsilon^{-1}$ this gives
\[
\E\norm{X-Z}^2=-2F'(\varepsilon^{-1})\le\varepsilon d\,F(\varepsilon^{-1})=\varepsilon d.
\]
For the second estimate, let $P$ be the Markov operator on $L^2(\mu)$ defined by
\[
Pf(x)=\E[f(Z)\mid X=x].
\]
For the coordinate functions $h_i(x)=x_i\in L^2(\mu)$,
\[
Ph_i(x)=\E[Z_i\mid X=x]=\bigl(T^\varepsilon_{\mu,\mu}(x)\bigr)_i.
\]
Thus
\[
\int\norm{x-T^\varepsilon_{\mu,\mu}(x)}^2\,\mu(dx)
=\sum_{i=1}^d\int |h_i(x)-Ph_i(x)|^2\,\mu(dx)
=\sum_{i=1}^d\norm{(I-P)h_i}_{L^2(\mu)}^2.
\]
Symmetry of the self-plan gives, for $f,g\in L^2(\mu)$,
\[
\langle f,Pg\rangle_{L^2(\mu)}=\E[f(X)g(Z)]=\E[f(Z)g(X)]=\langle Pf,g\rangle_{L^2(\mu)}.
\]
Thus $P$ is self-adjoint. Conditional Jensen's inequality gives
\[
|Pf(X)|^2=|\E[f(Z)\mid X]|^2
\le\E[|f(Z)|^2\mid X].
\]
Taking expectations and using that both marginals are $\mu$, we obtain
\[
\norm{Pf}_{L^2(\mu)}^2
=\E|Pf(X)|^2
\le\E\!\left[\E[|f(Z)|^2\mid X]\right]
=\E|f(Z)|^2
=\norm f_{L^2(\mu)}^2.
\]
Hence $\norm{Pf}_{L^2(\mu)}\le\norm f_{L^2(\mu)}$, and thus $P\preceq I$. We next show that $P\succeq0$. Since $d\alpha/d\mu$ is bounded, $f\in L^2(\mu)$ implies $f\in L^2(\alpha)$, and finiteness of $\alpha$ gives $f\in L^1(\alpha)$ by Cauchy--Schwarz. By the factorization~\eqref{eq:schrodinger-factorization},
\[
\langle f,Pf\rangle_{L^2(\mu)}
=\iint e^{-\norm{x-z}^2/(2\varepsilon)}f(x)f(z)\,\alpha(dx)\alpha(dz).
\]
Applying the calculation in \eqref{eq:self-fourier} to this double integral with $s=\varepsilon^{-1}$ and $f\alpha$ in place of $\alpha$ gives the identity below. The same absolute-integrability estimate justifies Fubini's theorem, with $(\int|f|\,d\alpha)^2<\infty$ in place of $\alpha(\R^d)^2$. Therefore, we have
\[
\langle f,Pf\rangle_{L^2(\mu)}
=\left(\frac{\varepsilon}{2\pi}\right)^{d/2}\int_{\R^d}
 e^{-\varepsilon\norm\xi^2/2}
 \left|\int e^{i\xi\cdot x}f(x)\alpha(dx)\right|^2d\xi\ge0.
\]
Hence $0\preceq P\preceq I$. Since $P$ and $I-P$ are commuting bounded positive self-adjoint operators, \cite[Theorem~9.3-1]{kreyszig} gives $P(I-P)\succeq0$. Consequently,
\[
(I-P)^2=I-P-P(I-P)\preceq I-P.
\]
For every $h\in L^2(\mu)$, self-adjointness of $I-P$ then gives
\[
\norm{(I-P)h}_{L^2(\mu)}^2
=\langle h,(I-P)^2h\rangle_{L^2(\mu)}
\le\langle h,(I-P)h\rangle_{L^2(\mu)}.
\]
Applying this inequality to the coordinate functions $h_i$ and summing over $i$ gives
\[
\int\norm{x-T^\varepsilon_{\mu,\mu}(x)}^2\,d\mu
=\sum_{i=1}^d\norm{(I-P)h_i}_{L^2(\mu)}^2
\le\int x\cdot(x-Px)\,d\mu.
\]
Since the self-plan is symmetric,
\[
\int x\cdot(x-Px)\,d\mu
=\E[X\cdot(X-Z)]
=\frac12\E\norm{X-Z}^2.
\]
Combining this with \eqref{eq:self} proves \eqref{eq:self-drift}.
\end{proof}

\section{Main proof and lemmas}

\subsection{Moment estimates and time-averaged transport bound}

Let
\[
m_p:=\int y\,p(dy),\qquad
M_p:=\int\norm{y-m_p}^2p(dy).
\]
Define the second moment about the target mean by
\[
M(t):=\int\norm{x-m_p}^2q_t(dx).
\]
\begin{lemma}\label{lem:moment}
At each time $t$, let $(X,Y)\sim\pi^\varepsilon_{q_t,p}$ and $(X,Z)\sim\pi^\varepsilon_{q_t,q_t}$. Then, for almost every $t$,
\begin{equation}
M'(t)=M_p-M(t)-\E_t\norm{X-Y}^2+\E_t\norm{X-Z}^2.
\label{eq:moment-id}
\end{equation}
Consequently, for every $t\ge0$,
\begin{equation}
\begin{aligned}
M(t)&\le e^{-t}M(0)+(M_p+\varepsilon d)(1-e^{-t})\\
&\le\max\{M(0),M_p+\varepsilon d\}.
\end{aligned}
\label{eq:moment-uniform}
\end{equation}
Moreover, for every $T>0$,
\begin{equation}
\frac1T\int_0^T\E_t\norm{X-Y}^2\,dt
\le M_p-\frac1T\int_0^TM(t)\,dt
+\varepsilon d-\frac{M(T)-M(0)}T.
\label{eq:cross-average-moment-bound}
\end{equation}
\end{lemma}

\noindent\emph{Remark.} To sample the variables in Lemma~\ref{lem:moment}, first draw $X\sim q_t$. Given its sampled value $x$, draw $Y$ and $Z$ independently from $\pi^\varepsilon_{q_t,p}(dy\mid x)$ and $\pi^\varepsilon_{q_t,q_t}(dz\mid x)$, respectively.

\begin{proof}
Since $q_t=(X_t)_\#q_0$,
\[
M(t)=\int\norm{x-m_p}^2q_t(dx)
=\int\norm{X_t(a)-m_p}^2q_0(da).
\]
Differentiating under the integral, justified by the support and velocity bounds on each finite time horizon, and using $\dot X_t(a)=V_t(X_t(a))$, we obtain
\[
\begin{aligned}
M'(t)
&=2\int(X_t(a)-m_p)\cdot\dot X_t(a)\,q_0(da)\\
&=2\int(X_t(a)-m_p)\cdot V_t(X_t(a))\,q_0(da)\\
&=2\int(x-m_p)\cdot V_t(x)\,q_t(dx)\\
&=2\int(x-m_p)\cdot\bigl(T^\varepsilon_{q_t,p}(x)-T^\varepsilon_{q_t,q_t}(x)\bigr)\,q_t(dx)\\
&=2\E_t\!\left[(X-m_p)\cdot\bigl(\E_t[Y-m_p\mid X]-\E_t[Z-m_p\mid X]\bigr)\right]\\
&=2\E_t[(X-m_p)\cdot(Y-m_p)]-2\E_t[(X-m_p)\cdot(Z-m_p)].
\end{aligned}
\]
Rewriting the inner products and using the marginals $X,Z\sim q_t$ and $Y\sim p$ gives
\[
\begin{aligned}
2\E_t[(X-m_p)\cdot(Y-m_p)]
&=\E_t\norm{X-m_p}^2+\E_t\norm{Y-m_p}^2-\E_t\norm{X-Y}^2\\
&=M(t)+M_p-\E_t\norm{X-Y}^2,\\
2\E_t[(X-m_p)\cdot(Z-m_p)]
&=\E_t\norm{X-m_p}^2+\E_t\norm{Z-m_p}^2-\E_t\norm{X-Z}^2\\
&=2M(t)-\E_t\norm{X-Z}^2.
\end{aligned}
\]
Subtracting gives \eqref{eq:moment-id}. Using \eqref{eq:self}, we obtain
\[
M'(t)\le-M(t)+M_p+\varepsilon d.
\]
Gronwall's inequality then gives \eqref{eq:moment-uniform}.
Finally, integrating \eqref{eq:moment-id} over $[0,T]$, dividing by $T$, and using \eqref{eq:self} gives \eqref{eq:cross-average-moment-bound}.
\end{proof}

To obtain the time-averaged transport bound, we also need to control the long-time average of the velocity divergence. The following lemma provides this control using the moment bound above and the entropy evolution identity.

\begin{lemma}\label{lem:div}
Let $(q_t)_{t\ge0}$ be the Sinkhorn flow under the assumptions of Theorem~\ref{thm:main}, with velocity field $V_t$. Then
\begin{equation}
\limsup_{T\to\infty}\frac1T\int_0^T\int\nabla_x\cdot V_t(x)q_t(dx)\,dt\le0.
\label{eq:average-divergence}
\end{equation}
\end{lemma}

\begin{proof}
Write $q_0(dx)=\rho_0(x)\,dx$. Since the assumptions do not guarantee that $\rho_0$ has finite entropy, we first truncate the initial law to justify the entropy evolution identity. For $n\ge1$, define
\[
A_n=\{x:\norm x\le n,\ 1/n\le\rho_0(x)\le n\}.
\]
Since $0<\rho_0<\infty$ for $q_0$-almost every $x$, $q_0(A_n)\uparrow1$. For $q_0(A_n)>0$, define the probability measures
\[
r_0^{(n)}=\frac{q_0|_{A_n}}{q_0(A_n)},\qquad r_t^{(n)}=(X_t)_\#r_0^{(n)},
\qquad r_t^{(n)}(dx)=\rho_t^{(n)}(x)\,dx.
\]
By construction, the truncated initial density $\rho_0^{(n)}$ has finite entropy, $\int\rho_0^{(n)}\log\rho_0^{(n)}\,dx<\infty$.
For every $T>0$, equation~\eqref{eq:velocity-divergence-bound} in the proof of Lemma~\ref{lem:density} gives a uniform bound on $\nabla_x\cdot V_t(x)$ for $t\in[0,T]$ and $x\in\supp q_t$. Equation~\eqref{eq:liouville-restricted} therefore implies that $\rho_t^{(n)}$ has finite entropy at every finite time. Taking logarithms of \eqref{eq:liouville-restricted} for $\rho_t^{(n)}$, integrating against $r_0^{(n)}$, and using $r_t^{(n)}=(X_t)_\#r_0^{(n)}$ gives
\begin{equation}
\int\rho_T^{(n)}\log\rho_T^{(n)}\,dx-\int\rho_0^{(n)}\log\rho_0^{(n)}\,dx
=-\int_0^T\int\nabla_x\cdot V_t\,dr_t^{(n)}\,dt.
\label{eq:core-entropy-identity}
\end{equation}
Moreover, $q_0(A_n)r_t^{(n)}=(X_t)_\#(q_0|_{A_n})\le (X_t)_\#q_0=q_t$, and therefore
\[
\int\norm{x-m_p}^2r_t^{(n)}(dx)\le\frac{M(t)}{q_0(A_n)}\le\frac{\max\{M(0),M_p+\varepsilon d\}}{q_0(A_n)}.
\]
The second inequality follows from \eqref{eq:moment-uniform} and is uniform in $t$. Comparing $r_t^{(n)}$ with the Gaussian density with mean $m_p$ and identity covariance
$\varphi_{m_p}(x)=(2\pi)^{-d/2}e^{-\norm{x-m_p}^2/2}$ gives
\[
0\le\KL(r_t^{(n)}\Vert\varphi_{m_p}(x)dx)
=\int\rho_t^{(n)}\log\rho_t^{(n)}\,dx
+\frac d2\log(2\pi)+\frac12\int\norm{x-m_p}^2r_t^{(n)}(dx).
\]
Consequently,
\begin{equation}
\int\rho_t^{(n)}\log\rho_t^{(n)}\,dx
\ge-\frac d2\log(2\pi)-\frac{\max\{M(0),M_p+\varepsilon d\}}{2q_0(A_n)}.
\label{eq:core-entropy-lower}
\end{equation} Dividing \eqref{eq:core-entropy-identity} by $T$, the finite initial entropy and \eqref{eq:core-entropy-lower} give
\begin{equation}
\limsup_{T\to\infty}\frac1T\int_0^T\int\nabla_x\cdot V_t\,dr_t^{(n)}\,dt\le0.
\label{eq:truncated-average-divergence}
\end{equation}

It remains to control the discarded mass. Put
\[
\lambda_t^{(n)}=(X_t)_\#(q_0|_{A_n^c}).
\]
Then $q_t=q_0(A_n)r_t^{(n)}+\lambda_t^{(n)}$ and $\lambda_t^{(n)}(\R^d)=1-q_0(A_n)$.
For $(X,Y)\sim\pi^\varepsilon_{q_t,p}$, we have the conditional second-moment identity
\[
\E_t[\norm{Y-m_p}^2\mid X=x]
=\tr\Cov_t(Y\mid X=x)+\norm{\E_t[Y\mid X=x]-m_p}^2.
\]
Together with \eqref{eq:jac-cov} and nonnegativity of the self-covariance trace, this yields
\[
\nabla_x\cdot V_t(x)
\le\varepsilon^{-1}\tr\Cov_t(Y\mid X=x)
\le\varepsilon^{-1}\E_t[\norm{Y-m_p}^2\mid X=x]
\le\frac{R_{p,c}^2}{\varepsilon},
\]
where $R_{p,c}:=\sup_{y\in\supp p}\norm{y-m_p}$ is finite by compact support.
Consequently,
\[
\int\nabla_x\cdot V_t\,d\lambda_t^{(n)}
\le\frac{R_{p,c}^2}{\varepsilon}(1-q_0(A_n))
\]
uniformly in $t$. Decomposing $q_t=q_0(A_n)r_t^{(n)}+\lambda_t^{(n)}$ and using \eqref{eq:truncated-average-divergence} yields
\[
\limsup_{T\to\infty}\frac1T\int_0^T\int\nabla_x\cdot V_t\,dq_t\,dt
\le\frac{R_{p,c}^2}{\varepsilon}(1-q_0(A_n)).
\]
Sending $n\to\infty$ proves the claim.
\end{proof}

\begin{proposition}[Time-averaged bound]\label{prop:cesaro-cross}
Along the flow of Theorem~\ref{thm:main}, let $(X,Y)\sim\pi^\varepsilon_{q_t,p}$ at each time $t$. Then
\begin{equation}
\limsup_{T\to\infty}\frac1T\int_0^T\E_t\norm{X-Y}^2\,dt\le\varepsilon d.
\label{eq:cross-cesaro}
\end{equation}
Consequently,
\[
\limsup_{T\to\infty}\frac1T\int_0^T W_2^2(q_t,p)\,dt\le\varepsilon d.
\]
\end{proposition}

\begin{proof}
By Lemma~\ref{lem:moment}, the averaged cross-cost satisfies \eqref{eq:cross-average-moment-bound}. The last term in that inequality tends to zero because $M$ is uniformly bounded by \eqref{eq:moment-uniform}.
Taking the limsup as $T\to\infty$ therefore gives
\[
\limsup_{T\to\infty}\frac1T\int_0^T\E_t\norm{X-Y}^2\,dt
\le M_p+\varepsilon d
-\liminf_{T\to\infty}\frac1T\int_0^T M(t)\,dt.
\]
Thus, to prove \eqref{eq:cross-cesaro}, it remains to show that
\[
\liminf_{T\to\infty}\frac1T\int_0^T M(t)\,dt\ge M_p.
\]
To obtain this lower bound, we relate $M(t)$ to the divergence of $V_t$. We have the conditional second-moment identities
\[
\begin{aligned}
\E_t[\norm{Y-m_p}^2\mid X=x]
&=\tr\Cov_t(Y\mid X=x)+\norm{\E_t[Y\mid X=x]-m_p}^2,\\
\E_t[\norm{Z-m_p}^2\mid X=x]
&=\tr\Cov_t(Z\mid X=x)+\norm{\E_t[Z\mid X=x]-m_p}^2.
\end{aligned}
\]
Since $X\sim q_t$, $Y\sim p$, and $Z\sim q_t$, integration gives
\[
\begin{aligned}
\int\E_t[\norm{Y-m_p}^2\mid X=x]\,q_t(dx)
&=\int\norm{y-m_p}^2p(dy)=M_p,\\
\int\E_t[\norm{Z-m_p}^2\mid X=x]\,q_t(dx)
&=\int\norm{z-m_p}^2q_t(dz)=M(t).
\end{aligned}
\]
Thus, using $\E_t[Y\mid X=x]=T^\varepsilon_{q_t,p}(x)$ and
$\E_t[Z\mid X=x]=T^\varepsilon_{q_t,q_t}(x)$, we obtain
\[
\begin{aligned}
\int\tr\Cov_t(Y\mid X=x)\,q_t(dx)
&=M_p-\int\norm{T^\varepsilon_{q_t,p}(x)-m_p}^2q_t(dx),\\
\int\tr\Cov_t(Z\mid X=x)\,q_t(dx)
&=M(t)-\int\norm{T^\varepsilon_{q_t,q_t}(x)-m_p}^2q_t(dx).
\end{aligned}
\]
Using $\nabla_x\cdot T=\tr(D_xT)$ and
$V_t=T^\varepsilon_{q_t,p}-T^\varepsilon_{q_t,q_t}$, equation~\eqref{eq:jac-cov}
and the two identities above give
\begin{equation}
\begin{aligned}
\varepsilon\int\nabla_x\cdot V_t\,dq_t
&=\int\tr\Cov_t(Y\mid X=x)\,q_t(dx)
-\int\tr\Cov_t(Z\mid X=x)\,q_t(dx)\\
&=\left[M_p-\int\norm{T^\varepsilon_{q_t,p}(x)-m_p}^2q_t(dx)\right]\\
&\quad-\left[M(t)-\int\norm{T^\varepsilon_{q_t,q_t}(x)-m_p}^2q_t(dx)\right]\\
&=M_p-M(t)+B(t),
\end{aligned}
\label{eq:divergence-moment}
\end{equation}
where
\[
B(t)=\int\bigl(\norm{T^\varepsilon_{q_t,q_t}(x)-m_p}^2-\norm{T^\varepsilon_{q_t,p}(x)-m_p}^2\bigr)q_t(dx).
\]
Since $V_t=T^\varepsilon_{q_t,p}-T^\varepsilon_{q_t,q_t}$,
\[
|B(t)|
\le\norm{V_t}_{L^2(q_t)}
\bigl(\norm{T^\varepsilon_{q_t,q_t}-m_p}_{L^2(q_t)}+\norm{T^\varepsilon_{q_t,p}-m_p}_{L^2(q_t)}\bigr).
\]
Conditional Jensen's inequality and the tower property give
\[
\begin{aligned}
\int\norm{T^\varepsilon_{q_t,q_t}-m_p}^2dq_t
&=\E_t\norm{\E_t[Z-m_p\mid X]}^2
\le\E_t\norm{Z-m_p}^2=M(t),\\
\int\norm{T^\varepsilon_{q_t,p}-m_p}^2dq_t
&=\E_t\norm{\E_t[Y-m_p\mid X]}^2
\le\E_t\norm{Y-m_p}^2=M_p.
\end{aligned}
\]
Thus
\[
|B(t)|\le\bigl(\sqrt{M(t)}+\sqrt{M_p}\bigr)\norm{V_t}_{L^2(q_t)}.
\]
Cauchy--Schwarz gives
\[
\begin{aligned}
\frac1T\int_0^T|B(t)|\,dt
&\le\frac1T\left(\int_0^T(\sqrt{M(t)}+\sqrt{M_p})^2\,dt\right)^{1/2}
\left(\int_0^T\norm{V_t}_{L^2(q_t)}^2\,dt\right)^{1/2}\\
&\le\frac{\sup_{t\ge0}(\sqrt{M(t)}+\sqrt{M_p})}{\sqrt T}
\left(\int_0^\infty\norm{V_t}_{L^2(q_t)}^2\,dt\right)^{1/2}
\longrightarrow0,
\end{aligned}
\]
since the supremum is finite by \eqref{eq:moment-uniform} and the last integral is finite by \eqref{eq:dissipation}.
Averaging and rearranging \eqref{eq:divergence-moment} gives
\[
\frac1T\int_0^T M(t)\,dt
=M_p+\frac1T\int_0^T B(t)\,dt
-\frac{\varepsilon}{T}\int_0^T\int\nabla_x\cdot V_t\,dq_t\,dt.
\]
The average of $B(t)$ tends to zero, and \eqref{eq:average-divergence} bounds the limsup of the average divergence by zero. Hence
\[
\liminf_{T\to\infty}\frac1T\int_0^TM(t)\,dt\ge M_p.
\]
Returning to \eqref{eq:cross-average-moment-bound} and taking the limsup now proves \eqref{eq:cross-cesaro}.
Finally, the cross-entropic plan is a coupling of $q_t$ and $p$, hence $W_2^2(q_t,p)\le \E_t\norm{X-Y}^2$.
\end{proof}

\subsection{A Wasserstein upper bound via the Sinkhorn divergence}\label{sec:coercivity}

We establish a static comparison without a target-moment term.

\begin{theorem}\label{thm:coercivity}
Let $\mu,\nu\in\Pt$. For every $\theta>1$,
\begin{equation}
2S_\varepsilon(\mu,\nu)
\ge\left(1-\frac1\theta\right)W_2^2(\mu,\nu)-\varepsilon d\log(2\theta).
\label{eq:coercivity}
\end{equation}
\end{theorem}

\begin{proof}
\noindent\textbf{Step 1. Construct a self-coupling.}\\
Let $(X,Y)\sim\pi^\varepsilon_{\mu,\nu}$ and put $C=\E\norm{X-Y}^2$.
 To bound the self-costs in terms of $C$, we use $Y$ together with a noisy observation of $X$ to construct a coupling of $\mu$ with itself. Fix $\sigma^2>0$, let $N\sim\mathcal N(0,\sigma^2I_d)$ be independent of $(X,Y)$, and set
\[
\widetilde X=X+N,\qquad U=(Y,\widetilde X).
\]
Conditionally on $U$, draw $X'$ independently with the same conditional law as $X$. Both marginals of $(X,X')$ equal $\mu$. Conditional independence gives the Markov chain $X-U-X'$, so the data-processing inequality for mutual information \cite{cover-thomas} gives
\[
I(X;X')\le I(X;U).
\]
Moreover, conditional independence and equality of the conditional laws give
\[
\begin{aligned}
\frac12\E\norm{X-X'}^2
&=\frac12\E\!\left[\E[\norm X^2\mid U]+\E[\norm{X'}^2\mid U]
-2\E[X\mid U]\cdot\E[X'\mid U]\right]\\
&=\E\!\left[\E[\norm X^2\mid U]-\norm{\E[X\mid U]}^2\right]\\
&=\E\norm{X-\E[X\mid U]}^2.
\end{aligned}
\]
Inserting this coupling in the variational definition of $\OT(\mu,\mu)$ gives
\begin{equation}
\OT(\mu,\mu)
\le\E\norm{X-\E[X\mid U]}^2+\varepsilon I(X;U).
\label{eq:posterior-bound}
\end{equation}

\noindent\textbf{Step 2. Bound $\E\norm{X-\E[X\mid U]}^2$.}\\
Since conditional expectation minimizes mean-square error,
\[
\E\norm{X-\E[X\mid U]}^2
\le\E\norm{X-g(U)}^2
\]
for every square-integrable estimator $g(U)$. Thus, it suffices to construct one estimator whose error can be bounded in terms of $C$ and $\sigma^2d$. Consider
\[
\widehat X_\lambda=(1-\lambda)Y+\lambda\widetilde X,\qquad \lambda\in\R.
\]
Since $N$ is independent of $(X,Y)$ and centered,
\[
\E\norm{X-\widehat X_\lambda}^2
=\E\norm{(1-\lambda)(X-Y)-\lambda N}^2
=(1-\lambda)^2C+\lambda^2\sigma^2d.
\]
Minimizing over $\lambda$ gives
\[
\lambda_*=\underset{\lambda\in\R}{\arg\min}\,
\bigl[(1-\lambda)^2C+\lambda^2\sigma^2d\bigr]
=\frac{C}{C+\sigma^2d},\qquad
\E\norm{X-\widehat X_{\lambda_*}}^2
=\frac{C\sigma^2d}{C+\sigma^2d}.
\]
Using this estimator in the opening inequality gives
\begin{equation*}
\E\norm{X-\E[X\mid Y,\widetilde X]}^2
\le\frac{C\sigma^2 d}{C+\sigma^2 d}.
\end{equation*}

\noindent\textbf{Step 3. Bound the mutual information $I(X;U)$ and combine the self-costs.}\\
By the chain rule for mutual information,
\[
I(X;Y,\widetilde X)=I(X;Y)+I(X;\widetilde X\mid Y).
\]
Writing $h$ for differential entropy, we obtain, for almost every $y$,
\[
\begin{aligned}
I(X;\widetilde X\mid Y=y)
&=h(\widetilde X\mid Y=y)-h(N)\\
&\overset{\text{\cite[Theorem~8.6.5, p.~254]{cover-thomas}}}{\le}
\frac12\log\det\!\left(2\pi e\,\Cov(\widetilde X\mid Y=y)\right)
-\frac d2\log(2\pi e\sigma^2)\\
&=\frac12\log\det\!\left(2\pi e\,[\Cov(X\mid Y=y)+\sigma^2I_d]\right)
-\frac d2\log(2\pi e\sigma^2)\\
&=\frac12\log\det\left(I_d+\frac{\Cov(X\mid Y=y)}{\sigma^2}\right)\\
&\le\frac d2\log\left(1+\frac{\tr\Cov(X\mid Y=y)}{\sigma^2 d}\right).
\end{aligned}
\]
The last inequality uses concavity of $\log$ applied to the eigenvalues.
Average in $y$ and use Jensen's inequality, we get
\[
\E\bigl[\tr\Cov(X\mid Y)\bigr]
=\E\norm{X-\E[X\mid Y]}^2
\le\E\norm{X-Y}^2=C.
\]
Hence
\[
I(X;\widetilde X\mid Y)
\le\frac d2\log\left(1+\frac{C}{\sigma^2 d}\right).
\]
The self-coupling bound \eqref{eq:posterior-bound} therefore yields
\[
\OT(\mu,\mu)\le\varepsilon I(X;Y)
+\frac{C\sigma^2 d}{C+\sigma^2 d}
+\frac{\varepsilon d}{2}\log\left(1+\frac{C}{\sigma^2 d}\right).
\]
Interchanging $X$ and $Y$ gives the same bound for $\OT(\nu,\nu)$. Subtracting the two self-cost bounds gives
\begin{equation}
\begin{aligned}
2S_\varepsilon(\mu,\nu)
&=2\OT(\mu,\nu)-\OT(\mu,\mu)-\OT(\nu,\nu)\\
&\ge C+2\varepsilon I(X;Y)
-2\left[\varepsilon I(X;Y)+\frac{C\sigma^2 d}{C+\sigma^2 d}
+\frac{\varepsilon d}{2}\log\left(1+\frac{C}{\sigma^2 d}\right)\right]\\
&=C-\frac{2C\sigma^2 d}{C+\sigma^2 d}
-\varepsilon d\log\left(1+\frac{C}{\sigma^2 d}\right).
\end{aligned}
\label{eq:coercivity-before-noise-choice}
\end{equation}
\noindent\textbf{Step 4. Choose the noise variance and conclude.}\\
We need to bound the nonlinear correction in \eqref{eq:coercivity-before-noise-choice} from above by an affine function of $C$.
Write that inequality as
\[
2S_\varepsilon(\mu,\nu)\ge C-H(C),
\]
where
\[
H(C)=\frac{2C\sigma^2 d}{C+\sigma^2 d}
+\varepsilon d\log\left(1+\frac{C}{\sigma^2 d}\right).
\]
A direct differentiation gives
\[
H''(C)=-\frac{4(\sigma^2 d)^2}{(C+\sigma^2 d)^3}
-\frac{\varepsilon d}{(C+\sigma^2 d)^2}<0,
\]
so $H$ is concave and every tangent line bounds it from above.
We seek a tangent with slope between zero and one, which we write as $1/\theta$ for a free parameter $\theta>1$: if
\[
H(C)\le\frac C\theta+b,
\]
with $b$ independent of $C$, then
\[
2S_\varepsilon(\mu,\nu)\ge\left(1-\frac1\theta\right)C-b,
\]
which retains a positive coefficient of $C$.
Choose $\sigma^2=\varepsilon\theta/(2\theta-1)$; at $C=\theta\varepsilon d$, this gives
\[
H(\theta\varepsilon d)=\varepsilon d\bigl(1+\log(2\theta)\bigr),
\qquad H'(\theta\varepsilon d)=\frac1\theta.
\]
The tangent line there therefore bounds $H$ from above:
\[
H(C)\le\varepsilon d\bigl(1+\log(2\theta)\bigr)
+\frac{C-\theta\varepsilon d}{\theta}
=\frac C\theta+\varepsilon d\log(2\theta).
\]
Substituting this bound for $H(C)$ into \eqref{eq:coercivity-before-noise-choice} gives
\[
2S_\varepsilon(\mu,\nu)
\ge\left(1-\frac1\theta\right)C-\varepsilon d\log(2\theta).
\]
Since $\theta>1$ and $W_2^2(\mu,\nu)\le C$, this proves \eqref{eq:coercivity}.
\end{proof}

\subsection{Limiting energy and asymptotic bound}

By \eqref{eq:edi}, the Sinkhorn divergence is nonincreasing along the flow; its nonnegativity therefore ensures that
\begin{equation}
S_\infty:=\lim_{t\to\infty}S_\varepsilon(q_t,p)\in\R.
\label{eq:energy-limit}
\end{equation}

\begin{lemma}\label{lem:energy}
Along the flow of Theorem~\ref{thm:main},
\begin{equation}
2S_\infty\le\left(1+\frac1{\sqrt2}\right)\varepsilon d.
\label{eq:limiting-energy}
\end{equation}
\end{lemma}

\begin{proof}
Fix $t\ge0$ and draw $(A,X)\sim\pi^\varepsilon_{q_t,q_t}$. Then
\[
\OT(q_t,q_t)=\frac12\E_t\norm{A-X}^2+\varepsilon I(A;X).
\]
Swapping $A$ and $X$ preserves the cost and the two marginals, so uniqueness of the optimizer gives $\mathcal L(A,X)=\mathcal L(X,A)$. Thus
\[
\E_t[A\mid X=x]=T^\varepsilon_{q_t,q_t}(x),
\qquad
\E_t[A-X\mid X]=T^\varepsilon_{q_t,q_t}(X)-X.
\]
Conditionally on $X$, draw $Y$ according to $\pi^\varepsilon_{q_t,p}(dy\mid X)$, independently of $A$. Then $(X,Y)\sim\pi^\varepsilon_{q_t,p}$. The Markov chain $A-X-Y$ gives, by data processing,
\[
I(A;Y)\le I(A;X).
\]
Using $\mathcal L(A,Y)$ as a coupling of $q_t$ and $p$ gives
\[
\begin{aligned}
\OT(q_t,p)
&\le\frac12\E_t\norm{A-Y}^2+\varepsilon I(A;Y)\\
&\le \OT(q_t,q_t)+\frac12\Bigl(\E_t\norm{A-Y}^2-\E_t\norm{A-X}^2\Bigr)\\
&=\OT(q_t,q_t)+\frac12\E_t\norm{X-Y}^2
+\E_t[(A-X)\cdot(X-Y)]\\
&=\OT(q_t,q_t)+\frac12\E_t\norm{X-Y}^2
+\E_t\bigl[\E_t[A-X\mid X]\cdot(X-Y)\bigr]\\
&=\OT(q_t,q_t)+\frac12\E_t\norm{X-Y}^2
-\E_t\bigl[(X-T^\varepsilon_{q_t,q_t}(X))\cdot(X-Y)\bigr].
\end{aligned}
\]
Repeating the argument with $(Y,X)\sim\pi^\varepsilon_{p,q_t}$ gives
\[
\OT(q_t,p)
\le\OT(p,p)+\frac12\E_t\norm{X-Y}^2
-\E_t\bigl[(Y-T^\varepsilon_{p,p}(Y))\cdot(Y-X)\bigr].
\]
Adding the two inequalities and subtracting $\OT(q_t,q_t)+\OT(p,p)$ gives
\begin{equation}
2S_\varepsilon(q_t,p)
\le\E_t\norm{X-Y}^2
-\E_t\bigl[(X-T^\varepsilon_{q_t,q_t}(X))\cdot(X-Y)\bigr]
-\E_t\bigl[(Y-T^\varepsilon_{p,p}(Y))\cdot(Y-X)\bigr].
\label{eq:transport-comparison}
\end{equation}

Since
\[
X-T^\varepsilon_{q_t,p}(X)=X-T^\varepsilon_{q_t,q_t}(X)-V_t(X),
\]
conditioning on $X$ gives
\[
\begin{aligned}
&-\E_t\bigl[(X-T^\varepsilon_{q_t,q_t}(X))\cdot(X-Y)\bigr]\\
&\quad=-\int\norm{x-T^\varepsilon_{q_t,q_t}(x)}^2q_t(dx)
+\int\bigl(x-T^\varepsilon_{q_t,q_t}(x)\bigr)\cdot V_t(x)q_t(dx)\\
&\quad\le\int\bigl(x-T^\varepsilon_{q_t,q_t}(x)\bigr)\cdot V_t(x)q_t(dx)\\
&\quad\le\left(\int\norm{x-T^\varepsilon_{q_t,q_t}(x)}^2q_t(dx)\right)^{1/2}\norm{V_t}_{L^2(q_t)}\\
&\quad\overset{\eqref{eq:self-drift}}{\le}\sqrt{\frac{\varepsilon d}{2}}\norm{V_t}_{L^2(q_t)}.
\end{aligned}
\]
Similarly, Cauchy--Schwarz and \eqref{eq:self-drift} give
\[
\begin{aligned}
-\E_t\bigl[(Y-T^\varepsilon_{p,p}(Y))\cdot(Y-X)\bigr]
&\le\left(\int\norm{y-T^\varepsilon_{p,p}(y)}^2p(dy)\right)^{1/2}
\left(\E_t\norm{X-Y}^2\right)^{1/2}\\
&\le\sqrt{\frac{\varepsilon d}{2}}\sqrt{\E_t\norm{X-Y}^2}.
\end{aligned}
\]
Substituting these bounds into \eqref{eq:transport-comparison} gives
\begin{equation}
2S_\varepsilon(q_t,p)
\le \E_t\norm{X-Y}^2+\sqrt{\frac{\varepsilon d}{2}}\norm{V_t}_{L^2(q_t)}
+\sqrt{\frac{\varepsilon d}{2}}\sqrt{\E_t\norm{X-Y}^2}.
\label{eq:energy-pointwise}
\end{equation}
Averaging \eqref{eq:energy-pointwise} over $[0,T]$ and applying Cauchy--Schwarz gives
\[
\begin{aligned}
\frac2T\int_0^T S_\varepsilon(q_t,p)\,dt
&\le\frac1T\int_0^T\E_t\norm{X-Y}^2\,dt\\
&\quad+\sqrt{\frac{\varepsilon d}{2T}}
\left(\int_0^\infty\norm{V_t}_{L^2(q_t)}^2\,dt\right)^{1/2}\\
&\quad+\sqrt{\frac{\varepsilon d}{2}}
\left(\frac1T\int_0^T\E_t\norm{X-Y}^2\,dt\right)^{1/2}.
\end{aligned}
\]
By \eqref{eq:energy-limit}, $\frac2T\int_0^T S_\varepsilon(q_t,p)\,dt\longrightarrow2S_\infty$ and
$\sqrt{\frac{\varepsilon d}{2T}}\left(\int_0^\infty\norm{V_t}_{L^2(q_t)}^2\,dt\right)^{1/2}\longrightarrow0$.
Taking the limsup and using \eqref{eq:cross-cesaro} gives
\[
2S_\infty
\le\varepsilon d+0+\sqrt{\frac{\varepsilon d}{2}}\sqrt{\varepsilon d}
=\left(1+\frac1{\sqrt2}\right)\varepsilon d,
\]
as claimed.
\end{proof}

\noindent\textbf{Proof of Theorem~\ref{thm:main}.}
Existence and uniqueness follow from Theorem~\ref{thm:flow}, and the time-averaged bound follows from Proposition~\ref{prop:cesaro-cross} since $W_2^2(q_t,p)\le\E_t\norm{X-Y}^2$ for $(X,Y)\sim\pi^\varepsilon_{q_t,p}$.

For the pointwise bound, fix any $\theta>1$.
Rearranging \eqref{eq:coercivity} at each finite time gives
\[
W_2^2(q_t,p)
\le\frac{\theta}{\theta-1}
\left[2S_\varepsilon(q_t,p)+\varepsilon d\log(2\theta)\right].
\]
Taking the limit superior as $t\to\infty$, using $S_\varepsilon(q_t,p)\to S_\infty$, and then applying Lemma~\ref{lem:energy}, we obtain
\begin{equation}
\begin{aligned}
\limsup_{t\to\infty}W_2^2(q_t,p)
&\le\frac{\theta}{\theta-1}
\left[2S_\infty+\varepsilon d\log(2\theta)\right]\\
&\le\frac{\theta}{\theta-1}
\left(1+\frac1{\sqrt2}+\log(2\theta)\right)\varepsilon d.
\end{aligned}
\label{eq:theta-asymptotic-bound}
\end{equation}
The minimizer over $\theta>1$ satisfies
\[
\theta_*-1=1+\frac1{\sqrt2}+\log(2\theta_*),
\qquad \theta_*\approx5.012111181619,
\]
and the resulting bound is
\[
\limsup_{t\to\infty}W_2^2(q_t,p)
\le\theta_*\,\varepsilon d
<5.012112\,\varepsilon d,
\]
which proves \eqref{eq:main}.
\hfill$\square$

\section{Conclusion}

We have established long-time Wasserstein error bounds for the Sinkhorn
Wasserstein gradient flow with a fixed regularization parameter
$\varepsilon>0$. For compactly supported initial and target laws,
with the initial law absolutely continuous, we obtain time-averaged
and asymptotic bounds of $\varepsilon d$ and $5.012112\,\varepsilon d$,
respectively, on the squared Wasserstein distance to the target.

Compact support provides a sufficient condition for global well-posedness.
We expect the same bounds to hold for $q_0,p\in\mathcal{P}_2(\mathbb{R}^d)$
with $q_0\ll\mathcal{L}^d$ whenever the Sinkhorn Wasserstein gradient flow
is globally well posed. While the present results provide quantitative
control of the long-time error, identifying conditions under which
the flow converges exactly to the target remains an important
direction for future research.

\end{document}